\documentclass{iig}

\usepackage[english]{babel}
\usepackage{amssymb,amsmath,eucal}
\usepackage[all,cmtip]{xy}
\usepackage{rotating}
\usepackage{tikz,relsize}
\newcommand{\fa}{\mathfrak{a}}
\newcommand{\fp}{\mathfrak{p}}

\newcommand{\mA}{\mathcal{A}}
\newcommand{\mD}{\mathcal{D}}

\newcommand{\Fun}{\mathbb{F}_1} %blocking sets

\newcommand{\Aut}{\mathrm{Aut}}
\newcommand{\id}{\mathbf{1}}

\newcommand{\Z}{\mathbb{Z}}
\newcommand{\GL}{\mathbf{GL}}
\newcommand{\F}{\mathbb{F}}
\newcommand{\uL}{\underline{\mathbb{L}}}

\def\doubleprod#1#2{\ooalign{$#1\prod$\cr$#1\coprod$\cr}}
\DeclareMathOperator*{\Rprod}{\mathpalette\doubleprod\relax}

\newcommand{\wis}[1]{{\text{\em \usefont{OT1}{cmtt}{m}{n} #1}}}
\def\Spec{\wis{Spec}}
\def\Sch{\wis{Sch}}

\newcommand{\mG}{\mathcal{G}} %grids

\newcommand{\K}{\mathbb{K}} %arcs, partial ovoids

\newcommand{\A}{\ensuremath{\mathbb{A}}}
\newcommand{\bP}{\ensuremath{\mathbb{P}}}
\newcommand{\bL}{\ensuremath{\mathbb{L}}}

\newcommand{\mO}{\ensuremath{\mathcal{O}}}

\newcommand{\PG}{\ensuremath{\mathbf{PG}}}
\newcommand{\PGL}{\ensuremath{\mathbf{PGL}}}

\newcommand{\AG}{\ensuremath{\mathbf{AG}}}

\newcommand{\mC}{\ensuremath{\mathcal{C}}}

\newcommand{\mE}{\mathcal{E}}
\newcommand{\mB}{\mathcal{B}}

\newcommand{\hP}{\wis{P}}
\newcommand{\B}{\mathbf{B}}

\begin{document}
\author[Thas]{Koen Thas}
\thanks{}
\address{Department of Mathematics \\
Ghent University\\
Krijgslaan 281, S25\\ 
B-9000 Ghent\\ 
Belgium}

\email{koen.thas@gmail.com}

\title{Counting points and acquiring flesh}
\date{}
\maketitle

\begin{abstract}
This set of notes is based on a lecture I gave at ``50 years of Finite Geometry | A conference on the occasion of Jef Thas's 70th birthday,'' in November 2014. It consists essentially of three parts: in a first part, I introduce some ideas which are based in the combinatorial theory underlying $\Fun$, the field with one element. In a second part, I describe,  in a nutshell,  the fundamental scheme theory over $\Fun$ which was designed by Deitmar. The last part focuses on zeta functions of Deitmar schemes, and also presents more recent work done in this area.
\end{abstract}


%\medskip
\setcounter{tocdepth}{1}
\tableofcontents

\medskip
\section{Introduction}

For a class of incidence geometries which are defined (for instance coordinatized) over  fields, it often makes
sense to consider the ``limit'' of these geometries when the number of field elements tends to $1$. As such, one ends up with
a guise of a ``field with one element, $\Fun$'' through taking limits of geometries. A general reference for $\Fun$ is the recent monograph \cite{KT-book}.

\subsection{Example: projective planes}
\label{explanes}

For instance, let the class of geometries be the 
classical projective planes $\PG(2,k)$ defined over commutative fields $k$. Then the number of points per line and the number of lines per point of such a plane is
\begin{equation}
\vert k \vert + 1, 
\end{equation}
so in the limit, the ``limit object'' should have $1 + 1$ points incident with every line. On the other hand, we want that the limit object {\em remains
an axiomatic projective plane}, so we still want it to have the following properties:
\begin{itemize}
\item[(i)]
any two distinct lines meet in precisely one point;
\item[(ii)]
any two distinct points are incident with precisely one line (the dual of (i));
\item[(iii)]
not all points are on one and the same line (to avoid degeneracy).
\end{itemize}

It is clear that such a limit projective plane (``defined over $\mathbb{F}_1$'') should be an ordinary triangle (as a graph).

\subsection{Example: generalized polygons}

Projective planes are, by definition,  generalized $3$-gons. Generalizing the situation to generalized $n$-gons, $n \geq 3$, a limit 
generalized $n$-gon becomes an ordinary $n$-gon (as a graph). The easiest way to see this is through the polygonal definition of generalized $n$-gons: if $\mE$ is the union of the set of points and the set of lines (which are assumed to be disjoint without loss of generality), then
one demands that: 
\begin{itemize}
\item[PD1]
there are no sub $m$-gons with $2 \leq m \leq n - 1$;
\item[PD2]
any two elements of $\mE$ are inside at least one $n$-gon, and
\item[PD3]
there exist $(n + 1)$-gons.
\end{itemize}

There is a constant $c \ne 0, 1$ such that any line is incident with $c + 1$ points \cite[1.5.3]{POL}, and as in the previous example, one lets $c$ go to $1$. So (PD3) cannot hold anymore.
In \cite{anal}, Tits defines a generalized $n$-gon over $\Fun$ to be an ordinary $n$-gon. (The fact that the number of lines incident with a point is also $1 + 1$, is explained at the end of \S\S \ref{Weylfunct}.)

\subsection{Example: Projective spaces of higher dimension}
\label{projhigh}

Generalizing the first example to higher dimensions, 
projective $n$-spaces over $\mathbb{F}_1$ should be sets $X$ of cardinality $n + 1$ endowed with the geometry of $2^X$: any subset (of cardinality $0 \leq r + 1 \leq n + 1$) is a subspace (of dimension $r$). 
In other words, projective $n$-spaces over $\Fun$ are {\em complete graphs on $n + 1$ vertices} with a natural subspace structure.
It is important to note that these spaces still satisfy the Veblen-Young axioms \cite{VY}, and that they are the only such incidence geometries with thin lines.

In the same vein, combinatorial affine $\Fun$-spaces consist of one single point and a number $m$ of one-point-lines through it; $m$ is the dimension of the space. We will come back to this definition in \S \ref{CT}.

In this paper, $\Aut(\cdot)$ denotes the automorphism group functor (from any category to the category of groups), and $\mathbf{S}_m$ denotes the symmetric group acting on $m$ letters.

\begin{proposition}[See, e.g., Cohn \cite{Cohn} and Tits \cite{anal}]
\label{propCT}
Let $n \in \mathbb{N} \cup \{-1\}$.
The combinatorial projective space $\PG(n,\mathbb{F}_1) = \PG(n,1)$ is the complete graph on $n + 1$ vertices endowed with the induced geometry of subsets,
and $\Aut(\PG(n,\mathbb{F}_1)) \cong \PGL_{n + 1}(\mathbb{F}_1) \cong \mathbf{S}_{n + 1}$. 
\end{proposition}

It is important to note that any $\PG(n,k)$ with $k$ a field contains (many) subgeometries isomorphic to $\PG(n,\Fun)$ as defined above; so the latter object is independent of $k$, and is the {\em common geometric substructure of all projective spaces of a fixed given dimension}:
\begin{equation}
\underline{\mA}: \{ \PG(n,k)\ \vert\ k \ \ \mbox{field} \} \longrightarrow \ \ \{ \PG(n,\Fun)\}.
\end{equation}

Further in this paper (in \S \ref{CT}), we will formally find the automorphism groups of $\Fun$-vector spaces through matrices, and these groups 
will perfectly agree with Proposition \ref{propCT}.

\subsection{Example: buildings}

The examples of the previous subsections can be generalized to all buildings $\mB$: in that case, the $\Fun$-copy is an apartment $\mA$. In the context of $\Fun$-geometry, apartments are often called {\em Weyl geometries}.
We refer to \cite{anal} and \cite{KT-combin} for details.

\subsection{Example: graphs}
\label{exgraphs}

Let $\Gamma$ be any graph, and see it as an incidence geometry with the additional property that any line/edge has precisely two distinct
points/vertices. (And let's assume for the sake of convenience that it has no loops.) Then over $\Fun$, nothing changes, and hence graphs 
are fixed points of the functor which sends incidence geometries to their $\Fun$-models.

\subsection{The functor $\underline{\mA}$}
\label{Weylfunct}

In \cite{KT-combin}, a functor $\underline{\mA}: \mG \mapsto \mG$ is described which associates to a natural class $\wis{B}$ of ``combinatorial $\Fun$-geometries'' $\mG$ its class $\wis{A}$ of ``$\Fun$-versions'' in much the same way as we have done here for the examples in \S\S \ref{explanes}|\S\S \ref{exgraphs}. These $\Fun$-versions can be obtained as fixed objects of $\underline{\mA}$ (which is called {\em Weyl functor} in {\em loc. cit.}).

The $\mathbb{F}_1$-functor $\underline{\mathcal{A}}$ should have several properties (with respect to the images); for the details, we refer to the chapter \cite{KT-combin}. Here, we isolate the following fundamental properties which will be useful for the present paper:

\begin{itemize}
\item[A1|]
all lines should have at most  $2$ different points;
\item[A2|]
an image should be a ``universal object,'' in the sense that it should be a subgeometry of any thick geometry of the same ``type'' (defined over {\em any} field, if at all defined over one) of at least the same rank; 
\item[A3|]
it should carry the same axiomatic structure (for example: $o \in \wis{A}$ and elements of $\underline{\mA}^{-1}(o)$ carry the same Buekenhout-Tits diagram);
\item[F|]
as $\wis{A}$ will be a subclass of the class of $\mathbb{F}_1$-geometries, it should consist precisely of the fixed elements of $\underline{\mathcal{A}}$.
 \end{itemize}

\br{\rm
We work up to point-line duality: that is why we are allowed to ask, without loss of generality, that lines have at most two points. We do {\em not} ask that they have {\em precisely}  two points, one motivation being e.g. (combinatorial) affine spaces over $\mathbb{F}_1$, in which any line has precisely one point. }
\er

In some sense, the number of lines through a point of an element $\Gamma$ of $\wis{A}$ should reflect the {\em rank} of the geometries in $\underline{\mathcal{A}}^{-1}(\Gamma)$. Think for example of the combinatorial affine and projective spaces over $\mathbb{F}_1$, and the ``Weyl geometries'' of buildings as described by Tits. This principle is a very important feature in the work of M\'{e}rida-Angulo and the author described in \S \ref{graphzeta}.

\section{Combinatorial theory}
\label{CT}

It is easy to see the symmetric group also directly as a limit with $\vert k \vert \longrightarrow 1$ of linear groups $\PG(n,k)$ (with the dimension fixed).
The number of elements in $\PG(n,k)$ (where $k = \F_q$ is assumed to be finite and $q$ is a prime power) is
\begin{equation}
\frac{(q^{n + 1} - 1)(q^{n + 1} - q)\cdots(q^{n + 1} - q^n)}{(q - 1)} = (q - 1)^nN(q)
\end{equation}
for some polynomial $N(X) \in \Z[X]$, and we have 
\begin{equation}
N(1) = (n + 1)! = \vert \mathbf{S}_{n +  1} \vert.\\
\end{equation}

\medskip
Now let $n,q \in \mathbb{N}$, and define $[n]_q = 1 + q + \cdots + q^{n - 1}$. (For $q$ a prime power, $[n]_q = \vert \PG(n,q)\vert$.) Put $[0]_q! = 1$, and define 
\begin{equation}
[n]_q! := [1]_q[2]_q\ldots [n]_q
\end{equation}
and
\begin{equation}
\left[\begin{array}{c}n \\ k\end{array}\right]_q = \frac{[n]_q!}{[k]_q![n - k]_q!}.
\end{equation}

If $q$ is a prime power, this is the number of $(k - 1)$-dimensional subspaces of $\PG(n - 1,q)$ ($=\vert \wis{Grass}(k,n)(\F_q)\vert$). The next proposition again gives sense to the limit 
situation of $q$ tending to $1$.

\begin{proposition}[See e.g. Cohn \cite{Cohn}]
The number of $k$-dimensional linear subspaces of $\PG(n,\mathbb{F}_1)$, with $k \leq n \in \mathbb{N}$, equals

\begin{equation}
\left[\begin{array}{c}n + 1 \\ k + 1\end{array}\right]_1 = \frac{n!}{(n - k)!k!} = \left[\begin{array}{c}n + 1 \\ k + 1\end{array}\right].
\end{equation}
\end{proposition}

\medskip
Many other enumerative formulas in Linear Algebra, Projective Geometry, etc. over finite fields $\F_q$ seem to keep meaningful interpretations if $q$
tends to $1$, and this phenomenon (the various interpretations) suggests a deeper theory in characteristic one.

Right now, we will have a look at some Linear Algebra features in characteristic $1$. Many of them are taken from Kapranov and Smirnov's \cite{KapranovUN}.

\subsection{A definition for $\Fun$}

One often depicts $\Fun$ as the set $\{0,1\}$ for which we only have the following 
operations:
\begin{equation}
0\cdot 1 = 0 = 0\cdot 0\ \ \mbox{and}\ \ 1\cdot 1 = 1. 
\end{equation}

This setting makes $\Fun$ sit in between the group $(\{ \id\},\cdot)$ and $\F_2$.
So in absolute Linear Algebra we are not allowed to have addition of vectors and we have to define everything in terms of scalar multiplication. 

The reason why this approach is natural, will become clear when we consider, e.g., linear automorphisms later in this section.

\subsection{Field extensions of $\Fun$}

For each $m \in \mathbb{N}^\times$ we define the {\em field extension}\index{field extension} $\mathbb{F}_{1^m}$\index{$\mathbb{F}_{1^m}$} of $\Fun$ of degree $m$ as the set
$\{0\} \cup \mu_m$, where $\mu_m$ is the (multiplicatively written) cyclic group of order $m$, and $0$ is an absorbing element for the extended
multiplication to $\{0\} \cup \mu_m$.

\subsection{Vector spaces over $\mathbb{F}_{1^{(n)}}$}

At the level of $\mathbb{F}_1$ we cannot make a distinction between 
affine spaces and vector spaces (as a torsor, nothing happens), so
in the vein of the previous section, a {\em vector/affine space}\index{vector space over $\F_{1^n}$}\index{affine space over $\F_{1^n}$} over $\F_{1^n}$, $n \in \mathbb{N}^\times$, is a triple $V = (\mathbf{0},X,\mu_n)$, where $\mathbf{0}$
is a distinguished point and $X$ a set, and
where $\mu_n$ acts freely on $X$. Each $\mu_n$-orbit corresponds to a direction. If $n = 1$, we get the notion considered in \S\S \ref{projhigh}.
If the dimension is countably infinite, $\mu_n$ may be replaced by $\mathbb{Z}, +$ (the infinite cyclic group). Another definition is needed when the dimension is larger.

\subsection{Basis}

A {\em basis}\index{basis} of the $d$-dimensional $\F_{1^n}$-vector space $V = (\mathbf{0},X,\mu_n)$ is a set of $d$ elements in $X$ which are 
two by two contained in different $\mu_n$-orbits (so it is a set of representatives of the $\mu_n$-action); here, formally, $X$ consists of $dn$ elements, and $\mu_n$ is the cyclic group with $n$ elements. (If $d$ is not finite one selects exactly one element in each $\mu_n$-orbit.)
If $n = 1$, we only have $d$ elements in $X$
(which expresses the fact that the $\mathbb{F}_1$-linear group indeed is the symmetric group) - as such we obtain the {\em absolute basis}\index{absolute!basis}.

Once a choice of a basis $\{b_i\ \vert\ i \in I\}$ has been made, any element $v$ of $V$ can be uniquely written as $b_j^{\alpha^u}$, for 
unique $j \in I$ and $\alpha^u \in \mu_n = \langle \alpha \rangle$. So we can also represent $v$ by a $d$-tuple with exactly one nonzero entry, namely $b_j^{\alpha^u}$ (in 
the $j$-th column).

\subsection{Dimension}

In the notation of above, the {\em dimension}\index{dimension of vector space} of $V$ is given by  $\mathrm{card}(V)/n = d$ (the number of $\mu_n$-orbits). 
%Note that if $q$ is a prime power, and $V$ is 
%an $\mathbb{F}_q$-vectorspace (say, of finite dimension), and $q \equiv 1 \mod{n}$, then 

\subsection{Field extension}

Let $V = (\mathbf{0},X,\mu_n)$ be a (not necessarily finite dimensional) $d$-space over $\F_{1^n}$, $n$ finite, so that $\vert X = X_V \vert = dn$. For any positive integral divisor $m$ of $n$, with $n = mr$,
$V$ can also be seen as a $dr$-space over $\F_{1^m}$. Note that there is a unique cyclic subgroup $\mu_m$ of $\mu_n$ of size $m$, so 
there is only one way to do it (since we have  to preserve the structure of $V$ in the process).

\subsection{Projective completion}

By definition, the {\em projective completion}\index{projective!completion} of an affine space $\AG(n,k)$, $n \in \mathbb{N}$ and $k$ a field, is the projective space $\PG(n,k)$ of the same dimension and defined over the same field, which one obtains by adding a hyperplane at infinity. 

We have seen how to perform projective completion over $\Fun$ through the following diagram:
\begin{equation}
\label{eqmotproj}
\PG(n,\Fun)\  =\  \AG(n,\Fun)\  +\  \PG(n - 1,\Fun).
\end{equation}

If one replaces $\Fun$ by an extension $\F_{1^m}$, the story is more complicate | see e.g. \cite{KT-combin,KT-proj}.

\subsection{Linear automorphisms}

A {\em linear automorphism}\index{linear!automorphism} $\alpha$ of an $\F_{1^n}$-vectorspace $V$ with basis $\{b_i\}$  is of the form 
\begin{equation}
\alpha(b_i) = b_{\sigma(i)}^{\beta_i} 
\end{equation}
for some power $\beta_i$ of the primitive $n$-th root of unity $\alpha$, and some permutation $\sigma \in \mathbf{S}_d$.
Then we have that 

\begin{equation}
\GL_d(\F_{1^n}) \cong \mu_n \wr \mathbf{S}_d.
\end{equation}

Elements of $\GL_d(\F_{1^n})$ can be written as $(d \times d)$-matrices with precisely one element of $\mu_n$ in each row or column (and
conversely, any such element determines an element of $\GL_d(\F_{1^n})$). In this setting, $\mathbf{S}_d$ is represented by $(d\times d)$-matrices with in each row and column exactly one $1$ | permutation matrices.

\br{\rm
Note that the underlying reason that rows and columns have only one nonzero element is that we do not have addition in our vector space.}
\er

\section{Deninger-Manin theory}

In a number of  works (\cite{Deninger1991}, \cite{Deninger1992} and \cite{Deninger1994}) on motives and regularized determinants, Deninger played with the possibility of translating Weil's proof of the Riemann Hypothesis for function fields of projective curves over finite fields $\F_q$ to the hypothetical curve $\overline{\Spec(\Z)}$. This idea also occurred, for instance, in Haran \cite{Haran}, and circulated in work of Smirnov \cite{Smirnov92} | see \cite{KT-Weil}. In \cite{Deninger1992}, Deninger gave a description of conditions on a certain category $\wis{M}$ of motives which might allow such a translation.

Let $\mC$ be a nonsingular absolutely irreducible projective algebraic curve over the finite field $\F_q$. Fix an algebraic closure $\overline{\F_q}$ of $\F_q$ and let $m \ne 0$ be a positive integer;  we have the following Lefschetz formula for the number $\vert \mC(\F_{q^m}) \vert$ of rational points over $\F_{q^m}$:

\begin{equation}
\vert \mC(\F_{q^m}) \vert = \sum_{\omega = 0}^2(-1)^{\omega}\mathrm{Tr}\Big(\mathrm{Fr}^m \Big| H^\omega(\mC)\Big) = 1 - \sum_{j = 0}^{2g}\lambda_j^m + q^f,
\end{equation}
where $\mathrm{Fr}$ is the Frobenius endomorphism acting on the \'{e}tale $\ell$-adic cohomology of $\mC$, the $\lambda_j$s are the eigenvalues of this action, and $g$ is the genus of the curve. We then have a motivic {\em weight decomposition}\index{weight!decomposition}
\begin{equation}
\zeta_{\mC}(s) = \prod_{\omega = 0}^2\zeta_{h^{\omega}(\mC)}(s)^{(-1)^{\omega - 1}} 
		= \frac{\prod_{j = 1}^{2g}(1 - \lambda_jq^{-s})}{(1 - q^{-s})(1 - q^{1 - s})} \nonumber \\
\end{equation}

%\begin{equation}
\begin{eqnarray} 
\label{compLefsch}
	 = \frac{\mbox{\textsc{Det}}\Bigl((s\cdot\id - q^{-s}\cdot\mathrm{Fr})\Bigl| H^1(\mC)\Bigr)}{\mbox{\textsc{Det}}\Bigl((s\cdot\id - q^{-s}\cdot\mathrm{Fr})\Bigl| H^0(\mC)\Bigr.\Bigr)\mbox{\textsc{Det}}\Bigl((s\cdot\id - q^{-s}\cdot\mathrm{Fr})\Bigl| H^2(\mC)\Bigr.\Bigr)}.
	\end{eqnarray}
	%\end{equation}
(Here the $\omega$-weight component is the zeta function of the pure weight $\omega$ motive $h^{\omega}(\mC)$ of $\mC$.)\\

The following analogous formula would hold in $\wis{M}$, where $\mC$ is replaced by the ``curve'' $\overline{\Spec(\Z)}$:

\begin{equation}
\label{eqDen}
\zeta_{\overline{\Spec(\Z)}}(s) = 2^{-1/2}\pi^{-s/2}\Gamma(\frac s2)\zeta(s)  
		=  \mathlarger{\frac{\Rprod_\rho\frac{s - \rho}{2\pi}}{\frac{s}{2\pi}\frac{s - 1}{2\pi}} } \overset{?}{=} \nonumber \\
\end{equation}
	\begin{eqnarray} 	
	 \frac{\mbox{\textsc{Det}}\Bigl(\frac 1{2\pi}(s\cdot\id - \varrho)\Bigl| H^1(\overline{\Spec(\Z)},*_{\mathrm{abs}})\Bigr.\Bigr)}{\mbox{\textsc{Det}}\Bigl(\frac 1{2\pi}(s\cdot\id -\varrho)\Bigl| H^0(\overline{\Spec(\Z)},*_{\mathrm{abs}})\Bigr.\Bigr)\mbox{\textsc{Det}}\Bigl(\frac 1{2\pi}(s\cdot\id - \varrho)\Bigl| H^2(\overline{\Spec(\Z)},*_{\mathrm{abs}})\Bigr.\Bigr)}. 
	\end{eqnarray}
	
	\bigskip
$\Big($The notation used in (\ref{eqDen}) is as follows:	
\begin{itemize}	
\item[$\ast$]
$\Rprod$ is the infinite {\em regularized product};
\item[$\ast$]
$\mbox{\textsc{Det}}$ denotes the {\em regularized determinant}\index{regularized determinant};
\item[$\ast$]
$\varrho$ is an ``absolute'' Frobenius endomorphism;
\item[$\ast$]
the $H^i(\overline{\Spec(\Z)},*_{\mathrm{abs}})$ are certain proposed cohomology groups, and
\item[$\ast$]
the $\rho$s run through the set of critical zeroes of the classical Riemann zeta.$\Bigr)$ 
\end{itemize}

Note that in the left-hand side of (\ref{eqDen}), we consider $\overline{\Spec(\Z)}$ instead of $\Spec(\Z)$, because we want to have a projective curve as in the expression for the motivic weight decomposition of $\mC$. This is 
why the factor 
\begin{equation}
2^{-1/2}\pi^{-s/2}\Gamma(\frac s2)
\end{equation}
occurs | it is the zeta-factor at infinity.

Conjecturally, in $\wis{M}$ there are motives $h^0$ (``the absolute point''), $h^1$ and $h^2$ (``the absolute Lefschetz motive'')\index{absolute!Lefschetz motive} with zeta functions
	\begin{equation}
	\label{eqzeta}
		\zeta_{h^w}(s) \ = \ \mbox{\textsc{Det}}\Bigl(\frac 1{2\pi}(s\cdot\id-\varrho)\Bigl| H^w(\overline{\Spec(\Z)},*_{\mathrm{abs}})\Bigr.\Bigr) 
	\end{equation}
for $w=0,1,2$. Deninger computed that $\zeta_{h^0}(s)=s/2\pi$ and $\zeta_{h^2}(s)=(s-1)/2\pi$. Manin proposed to interpret $h^0$ as $\Spec(\Fun)$ and $h^2$ as the affine line over $\Fun$, in \cite{Manin}. 

In \cite{Manin}, Manin then suggested to develop Algebraic Geometry over the field with one element, already in this specific context. 
So what {\em is} a scheme over $\Fun$?

\section{Deitmar schemes}

One of the first papers which systematically studied  a scheme theory over $\mathbb{F}_1$
was Deitmar's  \cite{Deitmarschemes2}, published in  2005. The study in \cite{Deitmarschemes2} is related to Kato's paper \cite{Kato}; see \S 5 and \S 9 of that paper.
By that time, Soul\'{e} had already published his fundamental $\mathbb{F}_1$-approach to varieties \cite{Soule}.

In $\Z$-scheme theory, a scheme $X$ is a locally ringed topological space which is locally isomorphic to affine schemes.
That is to say, $X$ is covered by opens $\{ U_i \ \vert \ i \in I \}$ such that the restriction of the structure sheaf $\mO_X$ to each $U_j$ is itself a locally
ringed space which is isomorphic to the spectrum of a commutative ring. 
 When aiming at an Algebraic Geometry over $\Fun$, one wants to have similar definitions at hand, but 
the commutative rings have to be replaced by  appropriate algebraic structures which reflect the $\Fun$-nature.

Several  attempts have been made to define schemes ``defined over $\mathbb{F}_1$,'' and often the approaches only differ in
subtle variations. 
We only need the most basic one, which is the ``monoidal scheme theory'' of Anton Deitmar \cite{Deitmarschemes1}. In this theory, 
the role of commutative rings over $\Fun$ is played by commutative monoids (with a zero).

\subsection{Rings over $\mathbb{F}_1$}

A {\em monoid}\index{monoid} is a set $A$ with a binary operation $\cdot: A \times A \longrightarrow A$ which is associative, and has an identity element (denoted $\id$).  Homomorphisms of monoids preserve units, and for a monoid $A$, $A^\times$\index{$A^\times$} will denote the group of invertible elements (so that if $A$ is a  group, $A^\times = A$).

In \cite{Deitmarschemes2}, Deitmar defines the category of rings over $\mathbb{F}_1$\index{$\Fun$-ring} to be the category of monoids (as thus ignoring additive structure), and the category of commutative $\Fun$-rings to be the category of commutative monoids. Usually, we will assume without further notice that an $\Fun$-ring $A$ als has a zero-element $0$ such that $0\cdot a = 0 = a \cdot 0$, $\forall a \in A$.

Below, all monoids will assumed to be abelian.

\subsection{Algebraic closure}

A monoid $A$ is {\em algebraically closed} if every equation of the form $x^n = a$ with $a \in A$ and $n \in \mathbb{N} \setminus \{0\}$ has $n$ solutions in $A$. Every monoid can be 
embedded into an algebraically closed monoid, and if $A$ is a group, then there exists a ``smallest'' such embedding which is called the {\em algebraic closure} of $A$.

The algebraic closure $\overline{\mathbb{F}_1}$\index{$\overline{\mathbb{F}_1}$} of $\mathbb{F}_1$ is the group $\mu_{\infty}$ of all complex roots of unity; it is isomorphic to $\mathbb{Q}/\Z$. Note that the multiplicative group $\overline{\mathbb{F}_p}^\times$ of the algebraic closure $\overline{\mathbb{F}_p}$ of 
the prime field $\mathbb{F}_p$ is isomorphic to the group of all complex roots of unity of order prime to $p$, so that the definition of $\overline{\mathbb{F}_1}$ might seem strange if compared with the finite field case. One can easily find ``meta-arithmetic'' arguments to deal with this matter | see \cite{KT-motive}.

\subsection{Localization}

Let $S$ be a submonoid of the monoid $A$. We define the monoid $S^{-1}A$\index{$S^{-1}A$}, the {\em localization}\index{localization} of $A$ {\em by} $S$, to be 
\begin{equation}
A \times S/\sim, 
\end{equation}
where the equivalence relation ``$\sim$'' is given by
\begin{equation}
(a,s) \sim (a',s') \ \ \mathrm{if\ and\ only\ if}\ \ s''s'a = s''sa'\ \ \mathrm{for\ some}\ \ s'' \in S.
\end{equation}

Multiplication in $S^{-1}A$ is componentwise, and one suggestively writes $\mathlarger{\frac{a}{s}}$ for the element in $S^{-1}A$ corresponding to $(a,s)$ (so $\mathlarger{\frac{a}{s}\cdot\frac{a'}{s'} = \frac{aa'}{ss'}}$).

\subsection{Ideal and spectrum}

If $C$ and $D$ are subsets of the monoid $A$, $CD$\index{$CD$} denotes the set of products $cd$, with $c \in C$ and $d \in D$. 

Recall that  a monoid is supposed to be abelian. If $C$ is a monoid, $\mathbb{Z}[C]$ denotes the corresponding ``monoidal ring'' | it is naturally defined similarly to a group ring.

An {\em ideal}\index{ideal} $\fa$ of a monoid $M$ is a
subset such that $M\fa \subseteq \fa$. For any ideal $\fa$ in $M$, $\mathbb{Z}[\fa]$ is an ideal in $\mathbb{Z}[M]$. Note that if $A$ and $B$
are monoids and $\alpha: A \longrightarrow B$ is a morphism, then $\alpha^{-1}(\fa)$ is an ideal in $A$ if $\fa$ is an ideal in $B$.

An ideal $\fp$  is called a {\em prime ideal}\index{prime!ideal}  if $S_{\fp} := M \setminus \fp$\index{$S_{\fp}$} is a monoid (that is, if $uv \in \fp$, then $u \in \fp$ or $v \in \fp$).
For any prime ideal $\fp$ in $M$, denote by  $M_{\fp} = S_{\fp}^{-1}M$
the {\em localization}\index{$M_{\fp}$}\index{localization} of $M$ {\em at} $\fp$. 

\begin{proposition}[Deitmar \cite{Deitmarschemes1}]
The natural map 
\begin{equation}
M \longrightarrow M_{\fp},\ m \longrightarrow \frac{m}{1}
\end{equation} 
with $\fp = M \setminus M^{\times}$ is an isomorphism.
\end{proposition}

Let $M$ be a monoid. The {\em spectrum}\index{spectrum} $\Spec(M)$\index{$\Spec(M)$} of $M$ is the set of prime ideals endowed with the obvious Zariski topology. 
Note that the spectrum cannot be empty since $M \setminus M^{\times}$ is a prime ideal.
The closed subsets are the empty set and all sets of the form

\begin{equation}
V(\fa) := \{ \fp \in \Spec(M) \vert \fa \subseteq \fp\},
\end{equation}\index{$V(\fa)$}
where $\fa$ is any ideal.
The point $\eta = \emptyset$ is contained in every nonempty open set and the point $M \setminus M^{\times}$ is closed and contained in every nonempty closed set.
Note also that for every $m \in M$ the set $V(m) :=  \{ \fp \in \Spec(M) \vert m \in \fp\}$ is closed (as $V(m) = V(Mm)$).

\begin{proposition}
$M \setminus M^{\times}$ is the unique maximal  ideal for any monoid $M$, so any such $M$ is a local $\mathbb{F}_1$-ring.
\end{proposition}

\subsection{Structure sheaf}

Let $A$ be a ring over $\mathbb{F}_1$. For any open set $U \subseteq \Spec(A)$, one defines $\mO_{\Spec(A)}(U) = \mO(U)$ to be the set of functions
(called {\em sections}\index{section})

\begin{equation}
s: U \longrightarrow \coprod_{\fp \in U}A_{\fp}
\end{equation}
for which $s(\fp) \in A_{\fp}$ for each $\fp \in U$, and such that there exists a neighborhood $V$ of $\fp$ in $U$, and elements $a, b \in A$, for which
$b \not\in \mathfrak{q}$ for every $\mathfrak{q} \in V$, and $s(\mathfrak{q}) = \mathlarger{\frac{a}{b}}$ in $A_{\mathfrak{q}}$. The map 

\begin{equation}
\mO_{\Spec(A)}: \Spec(A) \longrightarrow \mbox{monoids}: U \longrightarrow \mO(U)
\end{equation}
is the {\em structure sheaf}\index{structure sheaf} of $\Spec(A)$.

\begin{proposition}[Deitmar \cite{Deitmarschemes1}]
\label{globsec}
\begin{itemize}
\item[{\rm (i)}]
For each $\fp \in \Spec(A)$, the stalk $\mO_{\fp}$ of the structure sheaf is isomorphic to the localization of $A$ at $\fp$.
\item[{\rm (ii)}]
For global sections, we have $\Gamma(\Spec(A),\mO) := \mO(\Spec(A)) \cong A$.
\end{itemize}
\end{proposition}

\subsection{Monoidal spaces}

A {\em monoidal space}\index{monoidal space} is a topological space $X$ together with a sheaf of monoids $\mO_X$. 
Call a morphism of monoids $\beta: A \longrightarrow B$ {\em local}\index{local morphism} if $\beta^{-1}(B^{\times}) = A^{\times}$. A {\em morphism}\index{morphism!of monoidal spaces} between monoidal spaces $(X,\mO_X)$ and $(Y,\mO_Y)$ is defined naturally: it is a pair $(f,f^\#)$ with $f: X \longrightarrow Y$ a 
continuous function, and 
\begin{equation}
f^\#: \mO_Y \longrightarrow f_*\mO_X
\end{equation}
a morphism between sheaves of monoids on $Y$. (Here, $f_*\mO_X$\index{$f_*\mO_X$}, the {\em direct image sheaf}\index{direct image sheaf} on $Y$ induced by $f$, is defined by $f_*\mO_X(U) := \mO_X(f^{-1}(U))$ for all open $U \subseteq Y$.)
The morphism is {\em local}\index{local morphism} if each of the induced morphisms $f^\#_x: \mO_{Y,f(x)} \longrightarrow \mO_{X,x}$ is local.

\begin{proposition}[Deitmar \cite{Deitmarschemes1}]
\label{propmor}
\begin{itemize}
\item[{\rm (i)}]
If $A$ is any $\mathbb{F}_1$-ring, we have that the pair
$(\Spec(A),\mO_A)$ defines a monoidal space.
\item[{\rm (ii)}]
If $\alpha: A \longrightarrow B$ is a morphism of monoids, then $\alpha$ induces a morphism of monoidal spaces
\begin{equation}
(f,f^{\#}): \Spec(B) \longrightarrow \Spec(A),
\end{equation}
yielding a functorial bijection
\begin{equation}
\mathrm{Hom}(A,B) \cong \mathrm{Hom}_{\mathrm{loc}}(\Spec(B),\Spec(A)),
\end{equation}
where on the right-hand side we only consider local morphisms (hence the notation).
\end{itemize}
\end{proposition}

\subsection{Deitmar's $\mathbb{F}_1$-schemes}

As in the theory of rings, we have defined a structure sheaf $\mO_X$ on the topological space $X = \Spec(M)$, with $M$ a commutative monoid (with a zero). We define a {\em scheme} over $\mathbb{F}_1$ to be a topological space together with a sheaf of monoids, locally isomorphic to spectra of monoids in the above sense. The details are below.\\

{\bf Affine schemes}.\quad
An {\em affine scheme}\index{affine!$\Fun$-scheme} over $\mathbb{F}_1$ is a monoidal space which is isomorphic to $\Spec(A)$ for some monoid $A$. Such schemes are coined with the term {\em affine Deitmar schemes}\index{affine!Deitmar scheme} or also {\em $\mD$-schemes}\index{affine!$\mD$-scheme}
or {\em $\mD_0$-schemes}\index{affine!$\mD_{0}$-scheme}. (The ``$\mD$'' stands for ``Deitmar''; sometimes  the sub-index $0$ is added to stress that monoids have a zero in this context.)\\

{\bf General schemes}.\quad
A monoidal space $X$ is a {\em scheme}\index{$\Fun$-scheme} over $\mathbb{F}_1$ if for every point $x \in X$ there is an open neighborhood $U \subseteq X$ such that 
$(U,\mO_{X \vert U})$ is an affine scheme over $\mathbb{F}_1$. As in the affine case, we also speak of {\em $\mD$-schemes}\index{$\mD$-scheme} and {\em $\mD_0$-schemes}\index{$\mD_0$-scheme}.\\

A {\em morphism}\index{morphism!of $\mD_{(0)}$-schemes} of $\mD_{(0)}$-schemes is a local morphism of monoidal spaces. 
A point $\eta$ of a topological space is a {\em generic point}\index{generic point} if it is contained in every nonempty open set.

\begin{proposition}[Deitmar \cite{Deitmarschemes1}]
\begin{itemize}
\item[{\rm (i)}]
Any connected $\mD_0$-scheme has a unique gen\-eric point $\emptyset$, and morphisms between connected schemes map generic points to generic points.  
\item[{\rm (ii)}]
For an arbitrary $\mD_0$-scheme $X$, $\mathrm{Hom}(\Spec(\mathbb{F}_1),X)$ can be identified with the set of connected components of $X$.
\end{itemize}
\end{proposition}

\section{Acquiring flesh (1)}

Given an $\mathbb{F}_1$-ring $A$, {\em Deitmar base extension to $\mathbb{Z}$}\index{Deitmar base extension} is defined by

\begin{equation}
A \otimes \mathbb{Z} = A \otimes_{\mathbb{F}_1} \mathbb{Z} = \mathbb{Z}[A].
\end{equation}

Denote the functor of base extension by $\mathcal{F}(\cdot,\otimes_{\mathbb{F}_1}\mathbb{Z})$\index{$\mathcal{F}(\cdot,\otimes_{\mathbb{F}_1}\mathbb{Z})$}.

Conversely, we have a forgetful functor $\mathcal{F}$ which maps any (commutative) ring (with unit) to its (commutative) multiplicative monoid.

\begin{theorem}[Deitmar \cite{Deitmarschemes2}]
The functor  $\mathcal{F}(\cdot,\otimes_{\mathbb{F}_1}\mathbb{Z})$ is left adjoint to $\mathcal{F}$, that is, for every ring 
$R$ and every $\mathbb{F}_1$-ring $A$ we have that
\begin{equation}
\mathrm{Hom}_{\mathrm{Rings}}(A \otimes_{\mathbb{F}_1}\mathbb{Z},R) \cong \mathrm{Hom}_{\mathbb{F}_1}(A,\mathcal{F}(R)).
\end{equation}
\end{theorem}

One obtains a functor

\begin{equation}
X \longrightarrow X_{\mathbb{Z}}
\end{equation}
from $\mD_0$-schemes to $\mathbb{Z}$-schemes, thus extending the base change functor 
in the following way:
(a) write a scheme $X$ over $\mathbb{F}_1$ as a union of affine $\mD_0$-schemes, $X = \cup_i \Spec(A_i)$; (b) then map it to
$\cup_i\Spec(A_i \otimes_{\mathbb{F}_1} \mathbb{Z})$ (glued via the gluing maps from $X$).

Similarly to the general case, we say that the $\mD_0$-scheme $X$ is {\em of finite type}\index{finite!type $\mD_0$-scheme} if it has a finite covering by affine schemes $U_i = \Spec(A_i)$ such that the
$A_i$ are finitely generated. 

\begin{proposition}[Deitmar \cite{Deitmarschemes1}]
$X$ is of finite type over $\mathbb{F}_1$ if and only if $X_{\mathbb{Z}}$ is a $\mathbb{Z}$-scheme of finite type.
\end{proposition}

Conversely, one has a functor from monoids to rings, and it is left adjoint to the 
forgetful functor that sends a ring $(R,+,\times)$ to the multiplicative monoid $(R, \times)$. 
A scheme $X$ over $\mathbb{Z}$ can be written as a union of affine schemes 

\begin{equation}
X = \cup_{i}\Spec(A_i)
\end{equation}
for some set of rings $\{A_i\}$. Then map $X$ to $\cup_i\Spec(A_i,\times)$ (using the gluing maps from $X$) to obtain a functor from schemes over $\mathbb{Z}$ to schemes over $\mathbb{F}_1$ which extends the aforementioned forgetful functor. \\

The next theorem, which is due to Deitmar,  shows that integral $\mD_0$-schemes of finite type become {\em toric varieties}, once pulled to $\mathbb{C}$. (We won't define toric varieties here; we refer the reader to any standard text on these structures. Details can also be found in \cite{KT-motive}.)

\begin{theorem}[Deitmar \cite{Deitmartoric}]
Let $X$ be a connected integral $\mD_0$-scheme of finite type. Then every irreducible component of $X_\mathbb{C}$ is a toric variety. The components of 
$X_{\mathbb{C}}$ are mutually isomorphic as toric varieties.
\end{theorem}

Other scheme theories over $\Fun$ are known for which the base change functor to $\Z$ is ``more general.'' We refer to the monograph \cite{KT-book}, and the 
chapters therein, for a garden of such scheme theories.

\medskip
\section{Kurokawa theory}

One of the main tools to understand $\Fun$-schemes are their zeta functions. In this section, we define the Kurokawa zeta function, and we mention some interesting results taken from \cite{Kurozeta}. We first start with collecting some basic notions on arithmetic zeta functions.

\subsection{Arithmetic zeta functions}

Let $X$ be a scheme of finite type over $\mathbb{Z}$ | a {\em $\mathbb{Z}$-variety}\index{$\Z$-variety}. This means that $X$ has a finite covering of affine $\mathbb{Z}$-schemes $\Spec(A_i)$ with the $A_i$
finitely generated over $\mathbb{Z}$.
Recall that if $\widetilde{X}$ is an $k$-scheme, $k$ a field, a point $x \in \widetilde{X}$ is {\em $k$-rational}\index{$\F$-rational} if 
the natural morphism 
\begin{equation}
k \hookrightarrow k(x)
\end{equation}
is an isomorphism, with $k(x)$ the residue field of $x$ . (Note at this point that a homomorphism of fields $f \longrightarrow g$ is necessarily injective.)
A morphism
\begin{equation}
\Spec(L) \longrightarrow \widetilde{X},
\end{equation}
with $L/k$ a field extension,
is completely determined by the choice of a point $x \in \widetilde{X}$ (namely the image of $\Spec(L)$ in $\widetilde{X}$) and a field extension $L/k(x)$ through the natural $k$-embedding
\begin{equation}
k(x) \hookrightarrow L.
\end{equation}
Whence the set of $L$-rational points of $\widetilde{X}$ can be identified with 
\begin{equation}
\mathrm{Hom}(\Spec(L),\widetilde{X}).
\end{equation}
(If $\widetilde{X} \cong \Spec(A)$ is affine, $A$ being a commutative ring, one also has the identification with $\mathrm{Hom}(A,L)$.)

In the next proposition, a $k$-scheme $X \mapsto \Spec(k)$ is said to be {\em locally of finite type} (over $k$) if $X$ has a cover of open affine 
subschemes $\Spec(A_i)$, with all the $A_i$ finitely generated $k$-algebras.

\begin{proposition}[Closed and rational points]
\label{proparze}
\begin{itemize}
\item[{\rm (1)}]
Let $X$ be a $\Z$-scheme of finite type.
A point $x$ of $X$ is closed if and only if its residue field $k(x)$ is finite. (Note that $\vert k(x) \vert = \mathrm{dim}(\overline{\{x\}})$ as a closed subscheme.)
\item[{\rm (2)}]
Let $k = \overline{k}$ be algebraically closed, and let $\widetilde{X} \longrightarrow \Spec(k)$ be a $k$-scheme which is locally of finite type. Then a point $x$ is closed if and only if  it is $k$-rational.
\item[{\rm (3)}]
More generally, let $k$ be any field. Then a point $x$ of the $k$-scheme $\widetilde{X} \longrightarrow \Spec(k)$, which is again assumed to be locally of finite type, is closed if and only if the field extension $k(x)/k$ is finite. A closed point is $k$-rational if and only if $k(x) = k$.
\end{itemize}
\end{proposition}

%The philosophy is that, for instance when $\widetilde{X} \longrightarrow \Spec(\F_q)$ (with $q$ a prime power) is an $\F_q$-scheme, the $\F_q$-rational
%points can be found

Assume again that $X$ is an arithmetic scheme. 
Let $\overline{X}$ be the ``atomization''\index{atomization} of $X$; it is the set of closed points, equipped  with the discrete topology and the sheaf of fields $\{k(x)\ \vert\ x\}$. For $x \in \overline{X}$, let $N(x)$ be the cardinality of the finite field $k(x)$, that is, the {\em norm}\index{norm} of $x$.
Define the {\em arithmetic zeta function}\index{arithmetic zeta function} $\zeta_X(s)$ as 
\begin{equation}
\zeta_X(s) := \prod_{x \in \overline{X}}\frac{1}{1 - N(x)^{-s}}. 
\end{equation} 

\subsubsection*{Examples}

We mention four standard examples.

\begin{itemize}
\item[{\bf Dedekind}]
Let $X = \Spec(A)$, where $A$ is the ring of integers of a number field $\K$; then $\zeta_X(s)$ is the {\em Dedekind zeta function}\index{Dedekind zeta function} of $\K$.\\
\item[{\bf Riemann}]
Put $X = \Spec(\mathbb{Z})$; then $\zeta_X(s)$ becomes the classical {\em Riemann zeta function}\index{Riemann zeta function}.\\
\item[{\bf Affine sp.}]
With $\mathbb{A}^n(X)$ being the affine $n$-space over a scheme $X$, $n \in \mathbb{N}$, one has 
\begin{equation}
\zeta_{\mathbb{A}^n(X)} = \zeta_X(s - n).
\end{equation}
\item[{\bf Projective sp.}]
And with $\mathbb{P}^n(X)$ being the projective $n$-space over a scheme $X$, $n \in \mathbb{N}$, one has 
\begin{equation}
\zeta_{\mathbb{P}^n(X)} = \prod_{j = 0}^n\zeta_X(s - j).
\end{equation}
The latter can be obtained inductively by using the expression for the zeta function of affine spaces.
\end{itemize}

\subsection{Kurawa theory}

In \cite{Kurozeta}, Kurokawa says a scheme $X$ is of {\em $\Fun$-type}\index{scheme of $\Fun$-type} if its arithmetic zeta function $\zeta_X(s)$ can be expressed in the form
\begin{equation}
\zeta_X(s) = \prod_{k = 0}^n\zeta(s - k)^{a_k}
\end{equation}
with the $a_k$s in $\Z$. A very interesting result in \cite{Kurozeta} reads as follows:

\begin{theorem}[Kurokawa \cite{Kurozeta}]
Let $X$ be a $\Z$-scheme. The following are equivalent.
\begin{itemize}
\item[{\rm (i)}]
\begin{equation}
\zeta_X(s) = \prod_{k = 0}^n\zeta(s - k)^{a_k}
\end{equation}
with the $a_k$s in $\Z$.
\item[{\rm (ii)}]
For all primes $p$ we have
\begin{equation}
\zeta_{X\vert \F_p}(s) = \prod_{k = 0}^n(1 - p^{k - s})^{-a_k}
\end{equation}
with the $a_k$s in $\Z$.
\item[{\rm (iii)}]
There exists a polynomial $N_X(Y) = \sum_{k = 0}^na_kY^k$ such that
\begin{equation}
\#X(\F_{p^m}) = N_X(p^m) 
\end{equation}
for all finite fields $\F_{p^m}$.
\end{itemize}
\end{theorem}

Kurokawa defines the {\em $\Fun$-zeta function}\index{$\Fun$-zeta function} of a $\Z$-scheme $X$ which is defined over $\Fun$ as 
\begin{equation}
\zeta_{X\vert \Fun}(s) :=  \prod_{k = 0}^n(s - k)^{-a_k}
\end{equation}
with the $a_k$s as above. Define, again as above, the {\em Euler characteristic}\index{Euler characteristic}
\begin{equation}
\#X(\Fun) := \sum_{k = 0}^na_k.
\end{equation}

The connection between $\Fun$-zeta functions and arithmetic zeta functions is explained in the following theorem, taken from \cite{Kurozeta}.

\begin{theorem}[Kurokawa \cite{Kurozeta}]
Let $X$ be a $\Z$-scheme which is defined over $\Fun$. Then
\begin{equation}
\zeta_{X\vert \Fun}(s) =  \lim_{p \longrightarrow 1}\zeta_{X\vert \F_p}(s)(p - 1)^{\# X(\Fun)}.
\end{equation}
Here, $p$ is seen as a complex variable (so that the left hand term is the leading coefficient of the Laurent expansion of $\zeta_{X \vert \Fun}(s)$ around $p = 1$).
\end{theorem}

\subsubsection*{Examples}

For affine and projective spaces, we obtain the following zeta functions (over $\Z$, $\F_p$ and $\Fun$, with $n \in \mathbb{N}^\times$).
\begin{itemize}
\item[{\bf Affine sp.}]
$\zeta_{\A^n\vert \Z}(s) = \zeta(s - n);$\\

$\zeta_{\A^n\vert \F_p}(s) = {(1 - p^{n - s})}^{-1};$\\

$\zeta_{\A^n\vert \Fun}(s) = {(s - n)}^{-1},$\\

\item[{\bf Projective sp.}]
$\zeta_{\bP^n\vert \Z}(s)  =  \zeta(s)\zeta(s - 1)\cdots\zeta(s - n);$\\

$\zeta_{\bP^n\vert \F_p}(s) = \Bigl((1 - p^{-s})(1 - p^{1 - s})\cdots(1 - p^{n - s})\Bigr)^{-1};$\\

$\zeta_{\bP^n\vert \Fun}(s) = \Bigl(s(s - 1)\cdots(s - n)\Bigr)^{-1}.$

\end{itemize}

\medskip
\section{Graphs and zeta functions}
\label{graphzeta}

In this section we will introduce a new zeta function for (loose) graphs through $\Fun$-theory, following the work of \cite{MMA-KT}. 

In \cite{KT-Japan}, starting with a {\em loose graph} $\Gamma$, which is  a graph in which one also allows edges with $0$ or $1$ point,
I associated a Deitmar scheme $\wis{S}(\Gamma)$ to $\Gamma$ of which the closed points correspond to the vertices of $\Gamma$. 

Some features of $\wis{S}(\cdot)$:
\begin{itemize}
\item
Fundamental properties and invariants of the Deitmar scheme can  be obtained from the 
combinatorics of the loose graph, such as connectedness and the isomorphism class of the automorphism group. 
\item
A number of combinatorial $\Fun$-objects (such as combinatorial $\Fun$-pro\-ject\-ive space) are just loose graphs, and moreover,
the associated Deitmar schemes are precisely  the scheme versions in Deitmar's theory of these objects. 
\end{itemize}

Translation properties such as in the first item above,  were a main goal of the note \cite{KT-Japan}: trying to handle $\Fun$-scheme theoretic issues at the graph theoretic level (bearing in mind how some standard loose graphs should give rise to some standard Deitmar schemes).
After base extension, some basic properties of the ``real'' schemes might then be controlled by the loose graphs, etc.

The idea of the recent work \cite{MMA-KT} is now to associate a Deitmar scheme to a loose graph in a more natural way, 
and to show that, after having applied Deitmar's ($\cdot \otimes_{\Fun}\mathbb{Z}$)-functor, the obtained Grothendieck schemes are {\em defined over $\Fun$} in Kurokawa's sense. So they come with a Kurokawa zeta function, and that is the zeta function we associate to loose graphs.

As in \cite{MMA-KT}, we will call the modified functor  ``$\wis{F}$\index{$\wis{F}$}.''
It has to obey a tight set of rules, of which we mention a few  important ones:
\begin{itemize}
\item[{\bf Rule} \#1]
The loose graphs of the affine and projective space Deitmar schemes  should correspond to affine and projective space Deitmar schemes.
\item[{\bf Rule} \#2]
A vertex of degree $m$ should correspond locally to an affine space $\A^m$.
\item[{\bf Rule} \#3]
An edge without vertices should correspond to a multiplicative group.
\item[{\bf Rule} \#4]
``The loose graph is the map to gluing.''\\
\end{itemize}

\br{\rm
\begin{itemize}
\item
Because of Rule \#1, the pictures of Tits and Kapranov-Smirnov of affine and projective spaces over $\Fun$ are in agreement with the functor $\wis{F}$. (This was also the case for the functor $\wis{S}$.)
\item
In general, Rule \#2 does not hold for the functor $\wis{S}$. As we expressed at the end of the first section (in the discussion about the  
functor $\underline{\mA}$), this property is highly desirable though. 
\item
Rule \#3 implies that we have to work with a more general version of Deitmar schemes, since we allow expressions of type 
\begin{equation}
\Fun[X,Y]/(XY = 1) 
\end{equation}
(where the last equation generates a congruence on $\Fun[X,Y]$). In \cite{KT-Japan}, I only worked with Deitmar schemes, thus yielding a less natural approach to what the effect on deleting edges is on the corresponding schemes. By the way, $\Fun[X,Y]$ denotes the free abelian monoid generated multiplicatively by $X$ and $Y$, enriched with a zero.
\item
The last rule means that for any two vertices $u, v$ of a loose graph $\Gamma$, the intersection of the local affine spaces $\A_u$ and $\A_v$ which arise in $\wis{F}(\Gamma)$ as defined by Rule \#2, can be read from $\Gamma$. In general, this is a highly nontrivial game to play, as the examples and booby traps in \cite{MMA-KT} show.
\end{itemize}}
\er

For the  details, we refer the reader to \cite{MMA-KT}.

\subsection{The Grothendieck ring over $\Fun$}

Many of the formulas and calculations in \cite{MMA-KT} are expressed in the language of {\em Grothendieck rings}. 

\begin{definition}
The {\em Grothendieck ring of (Deitmar) schemes}\index{Grothendieck!ring! of $\Fun$-schemes} of finite type over $\Fun$, denoted as $K_0(\Sch_{\Fun})$\index{$K_0(\Sch_{\Fun})$}, is generated by the isomorphism classes of schemes ${X}$ of finite type over $\Fun$, $[X]_{_{\Fun}}$, with the relation
\begin{equation}
[X]_{_{\Fun}}= [X\setminus Y]_{_{\Fun}} + [Y]_{_{\Fun}} 
\end{equation}
for any closed subscheme $Y$ of $X$, and with the product structure given by
\begin{equation}
[X]_{_{\Fun}}\cdot[Y]_{_{\Fun}}= [X\times_{\Fun}Y]_{_{\Fun}}.
\end{equation}
\end{definition}
 
Denote by $\underline{\bL}=[\mathbb{A}^1_{\Fun}]_{_{\Fun}}$\index{$\underline{\bL}$} the class of the affine line over $\Fun$. Then the multiplicative group $\mathbb{G}_m$ 
satisfies
\begin{equation}
[\mathbb{G}_m]_{_{\Fun}}= \underline{\bL} - 1,
\end{equation} 
since it can be identified with the affine line minus one point.

If $X$ is a Deitmar scheme of finite type, and 
\begin{equation}
[X]_{\Fun}  \in\ \mathbb{Z}[\uL]\ \subset\ K_0(\Sch_{\Fun}), 
\end{equation}
then we say that $[X]_{\Fun} =: \wis{P}(X)$ is the {\em Grothendieck polynomial} of $X$.

\subsection{Affection principle}

Starting from a (finite) loose graph $\Gamma$, we denote the Deitmar scheme obtained by applying the functor $\wis{F}$ by $\wis{F}(\Gamma)$, as before. 

In \cite{MMA-KT} it is shown that $[\wis{F}(\Gamma)]_{\Fun} \in \mathbb{Z}[\uL]$.  Let $\hP(\Gamma)$ be the Grothendieck polynomial of $\wis{F}(\Gamma)$.
For each finite field $\mathbb{F}_q$, the number of $\F_q$-rational points of $\wis{F}(\Gamma) \otimes_{\Fun} \F_q$ is given by substituting the value $q$ for the indeterminate $\uL$ in $\hP(\Gamma)$ \cite{MMA-KT}. By {\bf Rule} \#3, locally each closed point of $\wis{F}(\Gamma)  \otimes_{\Fun} \F_q$ yields an affine space (of which the dimension is the degree of the point in the graph), 
so the total number of points can be expressed through the Inclusion-Exclusion principle. 

Consider a finite loose graph $\Gamma$, and let $\hP(\Gamma)$ be as above.  Taking any edge $uv$ which is not loose, we want to compare $\hP(\Gamma)$ and $\hP(\Gamma_{uv})$ in order to introduce a recursive procedure to simplify the loose graph (in that the number of cycles is reduced). Here, $\Gamma_{uv}$\index{$\Gamma_{uv}$} is the loose graph which one obtains when deleting the edge $uv$, while replacing it by two new loose edges, one through $u$ and one through $v$.

In this section, $\overline{\A}$ denotes the projective completion of the affine space $\A$. Also, if $\Gamma$ is a loose graph, $\bP(\Gamma)$ is the projective $\Fun$-space which is defined on the ambient graph of $\Gamma$ (i.e., the smallest graph in which $\Gamma$ is embedded). 

Calling $\mathrm{d}(\cdot,\cdot)$ the distance function in $\Gamma$ defined on $V \times V$, $V$ being the vertex set (so that, for example, $\mathrm{d}(s,t)$, with
$s$ and $t$ distinct vertices, is the number of edges in a shortest path from $s$ to $t$), it appears that one only needs to consider what happens 
in the vertex set
\begin{equation}
\B(u,1) \cup \B(v,1),
\end{equation}
where $\B(c,k) := \{ v \in V\ \vert\ \mathrm{d}(c,v) \leq k \}$. 

\medskip
\begin{theorem}[Affection Principle\index{Affection Principle} \cite{MMA-KT}]
\label{AP}
Let $\Gamma$ be a finite connected loose graph, let $xy$ be an edge on the vertices $x$ and $y$, and let $S$ be a subset of the vertex set. 
Let $k$ be any finite field, and consider the $k$-scheme $\wis{F}(\Gamma) \otimes_{\Fun}k$. Then $\cap_{s \in S}\A_s$ changes when one resolves 
the edge $xy$  only if $\cap_{s \in S}\A_s$ is contained in the projective subspace of $\bP(\Gamma) \otimes_{\Fun}k$ ``$k$-generated'' by $\B(x,1) \cup \B(y,1)$.
\end{theorem}

 In the next theorem, we will  use the notation $\bP(\B(u,1) \cup \B(v,1)) =: \bP_{u,v}$. If $\Delta$ is a loose graph, its {\em reduced}\index{reduced graph} version is the graph one obtains after deleting the loose edges.

\medskip
\begin{corollary}[Geometrical Affection Principle \cite{MMA-KT}]
\label{GAP}
Let $\Gamma$ be a finite connected loose graph, let $xy$ be an edge on the vertices $x$ and $y$, and let $k$ be any finite field.
The difference in the number of $k$-points of $\wis{F}(\Gamma) \otimes_{\Fun}k$ and $\wis{F}(\Gamma_{xy}) \otimes_{\Fun}k$ is 
\begin{equation}
\Bigl\vert \wis{F}(\Gamma_{\vert \bP_{x,y}}) \otimes_{\Fun}k \Bigr\vert_k\ -\ \Bigl\vert \wis{F}({\Gamma_{xy}}_{\vert \bP_{x,y}}) \otimes_{\Fun}k\Bigr\vert_k.
\end{equation}
In this expression, $\Gamma$ may be chosen to be reduced (but one is not allowed to reduce $\Gamma_{xy}$).
\end{corollary}

In terms of Grothendieck polynomials, we have the following theorem.

\begin{corollary}[Polynomial Affection principle \cite{MMA-KT}]
%\label{PAP}
Let $\Gamma$ be a finite connected loose graph, let $xy$ be an edge on the vertices $x$ and $y$, and let $k$ be any finite field.
Then in $K_0(\texttt{Sch}_k)$ we have
\begin{equation}
\hP(\Gamma) - \hP(\Gamma_{xy}) = \hP(\Gamma_{\vert \bP_{x,y}}) - \hP({\Gamma_{xy}}_{\vert \bP_{x,y}}).
\end{equation}
\end{corollary}

\subsection{Loose trees}
\label{loosetree}

Let $\Gamma$ be a loose tree. 

\begin{itemize}
\item Let $D$ be the set of degrees $\{d_1, \ldots, d_k \}$ of the vertex set $V(\Gamma)$ such that $1 < d_1 < d_2 < \ldots < d_k$.

\item Let us call $n_i$ the number of vertices of $\Gamma$ with degree $d_i$, $1\leq i \leq k$.

\item Put $\displaystyle I= \sum_{i=1}^k n_i - 1$.

\item Let $E$ be the number of vertices of $\Gamma$ with degree 1, that is the {\em end points}.
\end{itemize}

Then by \cite{MMA-KT}, the class of $\Gamma$ in $K_0(\Sch_{\Fun})$ is given by the  following map:

\begin{equation}
\label{formtree}
\begin{array}{c c c c}
\big[ \cdot \big]_{_{\Fun}} : & \{\mbox{Loose trees}\} & \longrightarrow  & K_0(\Sch_{\Fun})\\[1.5ex]
& \Gamma & \longrightarrow & \big[\Gamma\big]_{_{\Fun}} =  \displaystyle\sum_{i = 1}^k n_i\underline{\bL}^{d_i} - I\cdot\underline{\bL} + I + E.
\end{array}
\end{equation}

\subsection{Surgery}

Calculating Grothendieck polynomials of general loose trees is very complicated | see the many examples analyzed
in \cite{MMA-KT}. In {\em loc. cit.}, a procedure called ``surgery''\index{surgery} is introduced, which makes it possible to determine such 
polynomials by ``resolution of edges,''\index{resolution of edges} eventually reducing the calculation to the tree case, and this is a case which was resolved completely (cf. \S\S \ref{loosetree}). 

When $\Delta$ is a loose graph, and $e = xy$ is an edge with vertices $x$ and $y \ne x$, {\em resolving} $e$ means that one constructs the loose graph $\Delta_{xy} = \Delta_e$ as before, i.e., the adjacency between $x$ and $y$ is broken, and replaced by two new loose edges (one on $x$ and one on $y$). (Locally, the dimensions of the affine spaces at $x$ and $y$ remain the same, and the dimension of the ambient projective space of $\Delta$ increases.) 

In a nutshell, the following happens, starting from a finite loose graph $\Gamma = (V,E)$.
\bigskip
\begin{itemize}
\item[{\bf Spanning}]
Choose an arbitrary loose spanning tree $T$ (obviously defined)  in $\Gamma$.
\item[{\bf Resolution}]
Let $S$ be the set of egdes of $\Gamma$ not in $T$ which are not loose. Order $S = \{ e_1,\ldots,e_n\}$.
Now resolve all the edges in $S$, as follows:
\begin{equation}
\Gamma\ \longrightarrow\ \Gamma_{e_1}\ \longrightarrow\ {(\Gamma_{e_1})}_{e_2}\ \longrightarrow\ \cdots
\end{equation}
while keeping track of all the polynomial differences
\begin{equation}
\Big[ \hP(\Gamma_{e_1}) - \hP(\Gamma)\Big],\ \ \Big[ \hP({(\Gamma_{e_1})}_{e_2}) - \hP(\Gamma_{e_1})\Big],\ \ \ldots
\end{equation}
which one calculates using the Affection Principle.
\item[{\bf Reduction}]
Once one has resolved all the edges in $S$, we obtain a tree, and by \S\S \ref{loosetree} we know its Grothendieck polynomial. Now use the list of differences in the previous step to write down the Grothendieck polynomial of $\Gamma$.
\end{itemize}

In \cite{MMA-KT} it is shown that surgery is independent of the choice of the spanning tree, and of the order in which one chooses to resolve the edges.

\subsection{Resolving edges | two examples}

We now explain some examples.

\subsubsection{Example \#1}

We define $\Gamma(u,v;2)$,\index{$\Gamma(u,v;m)$} with $u, v$ symbols, to be the loose graph with adjacent
vertices $u, v$; $2$   common neighbors $v_1, v_2$ of $u$ and $v$ and no further incidences.

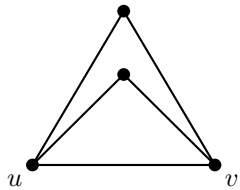
\begin{center}
\item
%\begin{figure}[h]
%\label{uvm}
  \begin{tikzpicture}[style=thick, scale=1.2]
\foreach \x in {-1,1}{
\fill (\x,0) circle (2pt);}

\fill (0,1) circle (2pt);
\fill (0,1.7) circle (2pt);
%\fill (0,2.4) circle (2pt);

\draw (-1,0) node[below left] {$u$} -- (1,0) node [below right] {$v$};
\draw (-1,0) -- (0,1);
\draw (0,1) -- (1,0);
\draw (-1,0) -- (0,1.7);
\draw (0,1.7) -- (1,0);
%\draw (-1,0) -- (0,2.4);
%\draw (0,2.4) -- (1,0);
%\draw[dotted] (-1,0) -- (-0.5,1.65);
%\draw[dotted] (1,0) -- (0.5,1.65);
%\draw[dotted] (-1,0) -- (-0.5,1.9);
%\draw[dotted] (1,0) -- (0.5,1.9);
\end{tikzpicture}
%\caption{The loose graph $\Gamma(u,v;m)$.}
%\end{figure}
\end{center}

For $k$ any field, the corresponding $k$-schemes consist of two affine $3$-spaces $\mathbb{A}_u$ and $\mathbb{A}_v$ and $2$ additional closed points in their spaces at infinity, of which the union covers all the points of the projective $3$-space $\bP(\Gamma(u,v;2))$ up to all points of the intersection $\gamma$ of their spaces at infinity (which is a projective line), except $2$ points in $\gamma$ in general position. So the Grothendieck polynomial is 
\begin{equation}
\label{eq6}
\sum_{i = 0}^{3}\uL^{i} - \Bigl((\sum_{i = 0}^{1}\uL^{i}) - 2\Bigr) = \uL^{3} + \uL^2 + 2.
\end{equation}

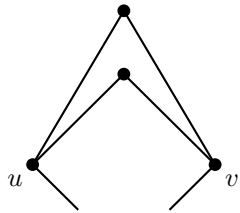
\begin{center}
%\begin{figure}[h]
%\label{res-uvm}
  \item
  \begin{tikzpicture}[style=thick, scale=1.2]
\foreach \x in {-1,1}{
\fill (\x,0) circle (2pt);}

\fill (0,1) circle (2pt);
\fill (0,1.7) circle (2pt);
%\fill (0,2.4) circle (2pt);

%\draw (-1,0) node[below left] {$u$} -- (1,0) node [below right] {$v$};
\draw (-1,0) node[below left] {$u$}  -- (0,1);
\draw (0,1) -- (1,0) node [below right] {$v$};
\draw (-1,0) -- (0,1.7);
\draw (0,1.7) -- (1,0);
%\draw (-1,0) -- (0,2.4);
%\draw (0,2.4) -- (1,0);
%\draw[dotted] (-1,0) -- (-0.5,1.65);
%\draw[dotted] (1,0) -- (0.5,1.65);
%\draw[dotted] (-1,0) -- (-0.5,1.9);
%\draw[dotted] (1,0) -- (0.5,1.9);

\draw (-1,0) -- (-0.5,-0.5);
\draw (1,0) -- (0.5,-0.5);
\end{tikzpicture}
%\caption{Resolution of $\Gamma(u,v;m)$ along the edge $uv$.}
%\end{figure}
\end{center}

Resolving $\Gamma(u,v;2)$ along $uv$, the $k$-schemes corresponding to $\Gamma(u,v;2)_{uv}$
consist of two disjoint affine $3$-spaces $\mathbb{A}_u$ and $\mathbb{A}_v$ (of which the planes at infinity intersect in the projective line generated by $v_1,v_2$) and $2$ additional mutually disjoint affine planes $\alpha_i$, $i = 1,2$,  in the projective $5$-space $\bP(\Gamma(u,v;2))$ such that for each $j$, $\alpha_j \cap \mathbb{A}_u \cong \alpha_j \cap \mathbb{A}_v$ is a projective line minus two points.

The Grothendieck polynomial is 
\begin{equation}
\label{eq8}
2\uL^{3} + 2\uL^2 - 4(\uL - 1).
\end{equation}

\subsubsection{Example \#2}

Starting from a triangle (as a graph), i.e., a combinatorial projective plane over $\Fun$, one deduces in a similar manner that its Grothendieck polynomial
is $\uL^2 + \uL + 1$. In general, the Grothendieck polynomial of a complete graph with $m + 1$ vertices, $m \ne 0$, is
\begin{equation}
\label{eqproj}
\uL^m + \uL^{m - 1} + \cdots + 1.
\end{equation}

The loose graph of an affine $\Fun$-space of dimension $m$ has as Grothendieck polynomial
\begin{equation}
\label{eqaff}
\uL^m.
\end{equation}

Both (\ref{eqproj}) and (\ref{eqaff}) are connected via the expression (\ref{eqmotproj}) in the Grothendieck ring.

\subsection{The zeta function}

We formally recall the next theorem (which was already mentioned implicitly), from \cite{MMA-KT}.

\begin{theorem}[\cite{MMA-KT}]
\label{MMAKTKuro}
For any loose graph $\Gamma$, the $\Z$-scheme $\chi := \wis{F}(\Gamma)\otimes_{\Fun}\Z$ is defined over $\Fun$ in Kurokawa's sense. 
\end{theorem}

Theorem \ref{MMAKTKuro} makes it possible to associate a (Kurokawa) zeta function to any loose graph, in the following way.

\begin{definition}[Zeta function for (loose) graphs]
Let $\Gamma$ be a loose graph, and let $\chi := \wis{F}(\Gamma) \otimes_{\Fun}\Z$. Let $P_{\chi}(X) = \sum_{i = 0}^ma_iX^i \in \Z[X]$ be as above.
We define the {\em $\Fun$-zeta function} of $\Gamma$ as:
\begin{equation}
\zeta^{\Fun}_{\Gamma}(t)\index{$\zeta^{\Fun}_{\Gamma}(t)$} := \displaystyle \prod_{k = 0}^m(t - k)^{-a_k}.\\
\end{equation}
\end{definition}

\medskip
\subsection{Example: trees}

Now let $\Gamma$ be a tree. We use the same notation as before, so that its class in the Grothendieck ring is given by

\begin{equation}
\big[\Gamma\big]_{_{\Fun}} =  \displaystyle\sum_{i = 1}^m n_i\underline{\bL}^{d_i} - I\cdot\underline{\bL} + I + E.
\end{equation}

The zeta function is thus given by
\begin{equation}
\zeta^{\Fun}_{\Gamma}(t)\ \ =\ \ \frac{(t - 1)^I}{t^{E + I}}\cdot\displaystyle \prod_{k = 1}^m(t - k)^{-n_k}.
\end{equation}

\medskip
\section{Acquiring flesh (2) | The Weyl functor depicted}

Sometimes, the functor $\underline{\mA}$ is artfully depicted by the following diagram, in which  Bacon's ``Study after Vel\'{a}zquez's portrait of Pope Innocent X'' \cite{Bacon}
 is compared to Vel\'{a}zquez's ``Portrait of Innocent X'' \cite{Vela} (Bacon's version being the $\Fun$-version of the original painting of Vel\'{a}zquez):

\begin{center}
\item
``Portrait of Innocent X''
\item
An oil on canvas ($114\mathrm{cm}\times 119\mathrm{cm}$) of the Spanish painter Diego Vel\'{a}zquez (1599-1660) dating from about 1650, depicting a portrait of Pope Innocent X. 
\end{center}

\begin{center}
\item
{\Large $\downarrow{{\small \underline{\mA}^{}}}$}
\end{center}

\begin{center}
\item
``Study after Vel\'{a}zquez's portrait of Pope Innocent X''
\item
An oil on canvas ($153\mathrm{cm}\times 118\mathrm{cm}$) of the Irish painter Francis Bacon (1909-1992) dating from 1953, showing a distorted version of Vel\'{a}zquez's  portrait of Pope Innocent X. 
\end{center}

\bigskip
At the conference, I showed that in a more modern setting, there is some analogy with the arrow 
\begin{center}
\item
{\large JAT\ \ {$\longrightarrow$}\ \ KT.}
\end{center}

\end{document}